# A Dynamic Holding Approach to Stabilizing a Bus Line Based on the Q-learning Algorithm with Multistage Look-ahead

Sheng-Xue He[*]

Business School, University of Shanghai for Science and Technology, No. 334, Jungong Road, Yangpu District, Shanghai 200093, China, lovellhe@usst.edu.cn

Jian-Jia He

Business School, University of Shanghai for Science and Technology, No. 334, Jungong Road, Yangpu District, Shanghai 200093, China, hejianjiayan@163.com

Shi-Dong Liang

Business School, University of Shanghai for Science and Technology, No. 334, Jungong Road, Yangpu District, Shanghai 200093, China, sdliang@hotmail.com

June Qiong Dong

School of Business, State University of New York at Oswego, 251 Rich Hall, 7060 State Route 104, New York, 13126-3599, USA, june.dong@oswego.edu

Peng-Cheng Yuan

Business School, University of Shanghai for Science and Technology, No. 334, Jungong Road, Yangpu District, Shanghai 200093, China, danis_cx@126.com

**Abstract:** The unreliable service and the unstable operation of a high frequency bus line are shown as bus bunching and the uneven distribution of headways along the bus line. Although many control strategies, such as the static and dynamic holding strategies, have been proposed to solve the above problems, many of them take on some oversimplified assumptions about the real bus line operation. So it is hard for them to continuously adapt to the evolving complex system. In view of this dynamic setting, we present an adaptive holding method which combines the classic approximate dynamic programming (ADP) with the multi-stage look-ahead mechanism. The holding time, that is the only control means used in this study, will be determined by estimating its impact on the operation stability of the bus line system in the remained observation period. The multi-stage look-ahead mechanism introduced into the classic Q-learning algorithm of the ADP model makes the algorithm get through its earlier unstable phase more quickly and easily. During the implementation of the new holding approach, the past experiences of holding operations can be cumulated effectively into an artificial neural network used to approximate the unavailable Q-factor. The use of a detailed simulation system in the new approach makes it possible to take into accounts most of the possible causes of instability. The numerical experiments show that the new holding approach can stabilize the system by producing evenly distributed headway and removing bus bunching thoroughly. Comparing with the terminal station holding strategies, the new method brings a more reliable bus line with shorter waiting times for passengers.

---

[*] Tel.: +086-18019495040.
E-mail address: lovellhe@usst.edu.cn

# 1. Introduction

If a bus is slowed down by some recurrent disturbances, it will likely encounter more passengers in the following stops and be delayed further; alternatively, a bus will speed up further if some reasons make it so at the beginning. In the first case, the bus slowed down will be caught up by the following bus; in the second case, the bus speeding up will catch up with its preceding bus. At last, buses will bunch up and cruise together. Newell and Potts (1964) first defined the above unpleasant phenomenon as "bus bunching". Bunching is an obvious symptom of the instability of bus line operation. Even if no bunching appears in a bus line system, the unevenly spaced buses or the uneven distribution of headways along the bus line manifests an unstable performance of the system. The negative influences of bus bunching or the uneven distribution of headways include the increased uncertainty of waiting time, the lengthened boarding time, the increased environmental pollution due to the inefficient operation of unevenly loaded buses, and the reduced attractiveness of transit system.

In order to resist bus bunching and improve the reliability of bus service, researchers have proposed many methods including stop-skipping strategies (Suh et al. 2002; Fu et al. 2003; Sun and Hickman 2005; Cortés et al. 2010; Liu et al. 2013), limited-boarding strategies (Osuna and Newell 1972; Newell 1974; Barnett 1974; Delgado et al., 2009, 2012), embedding-slack strategies (Daganzo 1997a, b; Daganzo 2009a; Xuan et al. 2011), static and dynamic holding strategies (Hickman 2001; Eberlein et al. 2001, Sun and Hickman 2008; Puong and Wilson 2008; Daganzo 2009a; Xuan et al. 2011; Bartholdi and Eisenstein 2012; Delgado et al. 2012; He 2015; Argote-Cabanero et al. 2015), speed adjustment strategies (Daganzo 2009b; Daganzo and Pilachowski 2009, 2011; He et al. 2019), transit signal priority strategies (Liu et al. 2003; Ling and Shalaby 2004; Estrada et al. 2016) and bus substitution strategy (Petit et al. 2018).

Many early analytical studies (Osuna and Newell 1972; Newell 1974; Barnett 1974; Ignall and Kolesar 1974; Hickman 2001; Zhao et al. 2006) mainly addressed the relatively simple systems. Generally only a single control point and either one or two buses are considered in these simple systems. With the advent of new technologies, some holding strategies have been proposed to take advantage of real-time information so as to reduce passengers' waiting times (Dessouky et al. 1999; Hichman 2001; Eberlein et al. 2001). Using real-time information, many headway-based holding strategies have been proposed to adaptively control the system (Daganzo 2009b; Bartholdi and Eisenstein 2012; Delgado et al. 2012; He 2015; Argote-Cabanero et al. 2015). But many of these strategies still need to simplify the reality by various assumptions to verify their effectiveness before being formulated as mathematical models. It is hard to integrate and deal with together in their formulations many important factors, such as signalized intersections, limited capacities of buses, random arrival rates of passengers, types of passenger behaviors and random travel times of road segments.

The headway-based adaptive holding strategies can use the real-time information to dynamically adjust the holding time to resist bus bunching (Daganzo 2009a; Xuan et al. 2011; He 2015). But once the criteria (or in another word the critical parameters) of the adjustment are determined, they will be remained during the observation period. On the one hand, it is hard to cumulate and learn from past experiences so as to make a better decision in the future. On the other hand, if some unusual events, e.g. surging arrivals appear at a stop or one bus breaks down, happen along the bus line, these holding strategies may be out of work.

In view of the complexity of a bus line system, system simulation is a popular tool used by many researches to deal with this issue (Koffman 1978; Turnquist and Bowman 1980; Abkowitz et al. 1986; Vandebona and Richardson 1986; Senevirante 1990; Adamski and Turnau 1998). But due to the limitation of the earlier theoretical development, the strategies are generally intuitive when simulation system is under consideration. It will be promising and beneficial to embed a detailed simulation system into some holding strategies which investigate in a more mathematical way various situations including the stochastic service, the real-time information, multiple control stops, etc (Hickman 2001; Eberlein et al. 2001, Sun and Hickman 2008; Puong and Wilson 2008; Delgado et al. 2012).

Since the daily travel pattern commonly repeats itself with small deviation day after day, the past experiences could be very helpful when used to choose a proper action to stabilize a high frequency bus line. When some unusual events happen during the observation period, a robust holding strategy should adapt itself properly and quickly to the new situation. Since the causes of bus line instability are very complicated including the passengers' arrival pattern at stops, the influence of other road traffic, the signal control schemes at intersections, departure frequency of buses, and so on, an effective holding strategy should take into account all the important factors at least in some implicit way. To realize the above expectations including making use of past experiences, dealing with unusual events, and covering most of the possible causes of instability, will make the holding strategy more robust and applicable. Or else, it is easy for a holding strategy to fall into the trap where it is effective only in some given situations.

The approximate dynamic programming (ADP) is a powerful tool that can be used to deal with systems with complex stochastic characteristics. Recently, Petit et al. (2018) have proposed a dynamic bus substitution strategy for bunching intervention based on ADP. One outstanding feature of ADP is that the past experiences can be learned by continuously refining the Q-factor function of the Q-learning algorithm of ADP during the observation period. In view of the above advantage of ADP, we will adopt it as our modeling basis. A reasonable ADP model for stabilizing a bus line needs some creative choices for its various components. This will be one important task of this paper. Holding a bus at a stop for a proper time period will be the only control means considered in this study. The proper holding time periods will be dynamically determined by solving the ADP model.

To deal with the unusual situations during the control period and later learn from the experience of dealing with these situations, the multi-stage look-ahead mechanism is introduced in the classic Q-learning algorithm to solve the ADP model. By looking ahead several succeeding stages, a proper holding time can be determined timely. The impact of some unusual events happened on the system stability can be timely and properly considered by the multi-stage look-ahead mechanism. This mechanism will also make the algorithm get through the earlier unstable phase more quickly.

One outstanding feature of the ADP model is that when a proper holding time is determined, the holding strategy based on ADP can take into account most of the possible causes of the instability, such as the passengers' arrival patterns and the signal control schemes at intersections. Firstly, a detailed simulation system can be adopted in the multi-stage look-ahead to facilitate the control decision making. Secondly, a detailed simulation system can be used as the deputy bus line system when the off-line version of the Q-learning algorithm of ADP is under consideration. Thirdly, the actual bus line system will be considered in the on-line version of the Q-learning algorithm. No matter which one, a simulation system or the actual system, is used, the various complicated causes of instability can be taken into account properly.

Two other features about the holding strategy to be investigated in this paper are worth mentioning here. The first one is that different from most of the holding strategies, the feasible boarding and alighting will be allowed during the whole dwell time in our new strategy. This will make the implementation of our strategy more acceptable. The second one is that though our strategy aims to stabilize the system (or equalize the uneven headways along the bus line) directly, the service level of bus line from the point of view of passengers will be improved as a byproduct.

The remainder of this paper is organized as follows. Section 2 introduces two indices for estimating the stability of bus line operation and some indices to evaluate the service level of bus line system. Section 3 formulates the ADP model by clarifying its main components. Section 4 presents the modified Q-learning algorithm with multi-stage look-ahead. Based on a detailed simulation system, Section 5 verifies our holding strategy with 5 bus lines with different scales. Finally, Section 6 provides concluding remarks.

## 2. Indices of System Stability and Service Level

Some basic assumptions and conceptions will be given in subsection 2.1. In subsection 2.2, we will introduce two indices to evaluate the operation stability of a bus line system over the whole observation period. In subsection 2.3, the commonly used service level indices of bus line systems, such as waiting time and running time of passengers, will be clarified from the point of view of passengers.

### 2.1. Basic Assumptions

To simplify the following analysis, we will use a circle bus line on which a given number of buses are cruising as our subject of study. Assume that boarding and alighting will separately take place through two different doors for all the buses under consideration. If the number of passengers in a bus reaches the maximal passenger capacity, boarding will be prohibited. But if there is any remained capacity, boarding and alighting are always allowed.

Holding a bus at a stop for a while is the only control means used in this study. Its starting and ending times are explained as follows. Assume that a bus $b$ arrives at a stop $e$ at time point $t$. The passengers arrived at $e$ before $t$ need to take time $t_1$ to board this bus. The passengers in this bus with $e$ as their destination will spend time $t_2$ to get off this bus. We assume that the possible holding operation on bus $b$ will start at the time point $t + \max\{t_1, t_2\}$. If this bus will be held at stop $e$ for a time interval $t_3$, the departure time of this bus from $e$ will be $t + \max\{t_1, t_2\} + t_3$. $t_3$ is the holding time that will be determined dynamically to stabilize the system. The detailed explanation of holding operation will be given in subsection 3.1.

### 2.2. Stability Indices

For a high frequency bus line, it is hard for buses to cruise along the bus line strictly following a fixed time table. The uneven headways of buses to their nearest preceding buses reflect the operation instability of the bus line. In order to measure the intensity of this operation stability, some proper indices need to be defined at first.

Firstly, we define the Dynamic Circle Headway (DCH) at time $t$ for a circle bus line as

follows:

$$H(t) = (\sum_{b \in B} h_b(t))/n_B, \tag{1}$$

where $h_b(t)$ is the time headway between bus $b$ and its nearest preceding bus. In another word, $h_b(t)$ is the expected time for bus $b$ to reach the current position of its preceding bus from its current position. $B$ is the set of all cruising buses and $n_B$ is the total number of buses in $B$. The value of $H(t)$ can be viewed as the expected instantaneous headway at time $t$ for a bus line with evenly dispersed buses.

With DCH, we can define the pseudo standard deviation of instantaneous headways at time $t$ as follows:

$$\sigma_H(t) = \sqrt{[\sum_{b \in B}(h_b(t) - H(t))^2]/n_B}. \tag{2}$$

If buses are evenly dispersed in terms of time along the bus line at time $t$, the value of $\sigma_H(t)$ will be small; or else, the value will be big.

During the observation period, we view the departure time of a bus from a bus stop as an estimating time point and group all such departure time points into a set $\bar{T}$. To simplify the following analysis, a departure time will be identified by the combination of its value, its corresponding bus and bus stop. So even two departure times have the same value, they are different departure times due to the fact that their corresponding buses or bus stops are different. We can calculate the average value of $\sigma_H(t)$ over $\bar{T}$ as follows:

$$\bar{c}_H = (\sum_{t \in \bar{T}} \sigma_H(t))/n_{\bar{T}}, \tag{3}$$

where $n_{\bar{T}}$ is the cardinality of $\bar{T}$. The corresponding standard deviation of $\bar{c}_H$ is

$$\sigma_{\bar{c}} = \sqrt{[\sum_{t \in \bar{T}}(\sigma_H(t) - \bar{c}_H)^2]/(n_{\bar{T}} - 1)}. \tag{4}$$

From the above definitions, we can see that $\bar{c}_H$ and $\sigma_{\bar{c}}$ have the same measurement unit as time. $\bar{c}_H$ stands for the deviation of headways from DCH over the whole observation period. $\sigma_{\bar{c}}$, i.e. the standard deviation of $\bar{c}_H$, shows the reliability of $\bar{c}_H$ as an expected estimate. $\bar{c}_H$ and $\sigma_{\bar{c}}$ are the stability indices to be used to evaluate the bus line's operation. From now on, $\bar{c}_H$ and $\sigma_{\bar{c}}$ will be referred to as the First Stability Index (**FSI**) and the Second Stability Index (**SSI**), respectively. Since $\bar{c}_H$ and $\sigma_{\bar{c}}$ are defined as average values over the cruising buses and the set $\bar{T}$ of departure times, they intrinsically possess the potential to indicate the operational stability of any bus line system as a whole.

**2.3. Indices of Service Level**

The waiting and riding times of passengers are the most important measures commonly used as the service level indices of a bus line. To properly measure the operational performance of a high frequency bus line at the end of the observation period, we need to divide the passengers into three groups and consider their time attributes separately.

There are three types of passengers when the observation reaches its end. The first type denoted by $P_1$ includes the passengers who have finished their travels by bus, i.e. having alighted at their destination stops during the observation period. The second type denoted by $P_2$ includes the passengers who are still in the buses at the end of the observation period. The third type denoted by $P_3$ consists of all the passengers still waiting at their original stops when the observation ends.

A passenger $p$ may have three time attributes including the out-of-vehicle waiting time at his/her original stop, the riding time and the total travel time. For the passengers in $P_i$ ($i = 1, 2, 3$), we will use $t_{P_i}^W$, $t_{P_i}^R$ and $t_{P_i}^{Tr}$ to denote the average waiting time, the average riding time and the average travel time over $P_i$. Obviously, $t_{P_i}^{Tr} = t_{P_i}^W + t_{P_i}^R$. $\sigma_{P_i}^W$, $\sigma_{P_i}^R$ and $\sigma_{P_i}^{Tr}$ denote the standard deviations of $t_{P_i}^W$, $t_{P_i}^R$ and $t_{P_i}^{Tr}$ over set $P_i$, respectively. Note that for a given type of passengers, some quantities mentioned above may be unavailable at the end of the observation. For example, $t_{P_i}^R$ and $t_{P_i}^{Tr}$ are unavailable for $P_3$. The time attributes and their corresponding deviations will be used as the service level indices later in the numerical analysis.

To correctly evaluate the service level of a bus line, the states of these three types of passengers should be evaluated separately and be considered together to estimate the service level of the whole system. For example, when only the first type is considered for a bus line with 10 to 15 cruising buses, the influence of the second type may be great, especially when the observation period is short.

The standard deviations of the above time attributes are required to distinguish the service levels indicated by the same time attributes. For example, if two operational results obtained from the same bus line have the same average waiting time, the one with the smaller standard deviation embodies a higher service level than the other.

**3. ADP Model**

The holding approach about to be presented is based on Approximate Dynamic Programming (ADP). Generally speaking, an ADP model consists of five main components. They are the control action, the information factor, the state variable, the state value function (or Q-factor) and the decision-making policy. We will define and explain these components one by one in the context of the bus line operation.

### 3.1. Set of Control Actions

In this study, we only use holding a bus at its current stop for a specified time interval as the unique control means of stabilizing a bus line. The holding operation starts at the current stop after both the boarding process for the passengers who have arrived at the stop before the bus arrives and the alighting process for the onboard passengers end. The operation of holding a bus at a stop for a specified time interval is defined as the control action. Note that the seemingly different concepts consisting of holding time, slack, control action and outside interference in the existing literatures generally have the same meaning.

A very short time interval $\tau$, e.g. 2 seconds, will be used as the measurement unit to measure the holding time. Assume that only the time interval with the form of $n\tau$ ($n$ is a non-negative integer) can be adopted as a feasible holding time in this study. We denote an action with holding time $n\tau$ by $a_{\to n}$. All the available actions constitute the action set $A \equiv \{a_{\to 0}, a_{\to 1}, \cdots, a_{\to n_A}\}$. It is possible and convenient to define a specified action set $A_e$ for a given bus stop $e$ such that $A_e \subseteq A$. Obviously, if there is only one action $a_{\to 0}$ in $A_e$, it means there is no holding at $e$.

To adopt an action set $A$ with finite elements as defined above can effectively avoid the curse of dimensionality of ADP with infinite actions. In the real world only the holding operation with reasonable duration can be carried out and accepted by most of passengers. Generally speaking, to set the time interval $\tau$ to 1 or 2 seconds is applicable in the real world.

The starting time to hold a bus at its current stop indicates a stage of the system operation. Note that a stage should be identified by the triples of the starting time of the holding operation, the bus to be held, and the bus stop where holding happens. The set of stages included in the observation period is denoted by $M$. Note that if two or more stages have the same starting time, the order to deal with them during the control period will be arbitrarily determined. A typical stage is indicated by $m$. The control action at stage $m$ is denoted by $a(m)$ or $a_m$. The cumulative intensity $a_\Sigma$ of the outside interference to the bus line is defined by the sum $\sum_{m \in M} a(m)$ of holding times. The average interference intensity over $M$ can be obtained by $\bar{a} = a_\Sigma / n_M$. The standard deviation of $\bar{a}$ is given by $\sigma_{\bar{a}} = \sqrt{[\sum_{m \in M}(a(m) - \bar{a})^2]/(n_M - 1)}$. Here $n_M$ is the cardinality of $M$. To resist bus bunching and improve the system stability with holding as the only control action, we hope the average interference intensity and its deviation as small as possible.

Since a bus will depart from a bus stop after the possible holding action has been carried out, the elements in set $\bar{T}$ of all the departure times correspond to the ones in $M$. In another word, there is a bijective mapping between $M$ and $\bar{T}$.

### 3.2. Information Factors

In terms of ADP, an information factor stands for a stochastic or uncertain factor influencing the system evolution. During the observation period, the uncertain factors will be fixed gradually and show their final influences with the evolution of the bus line system. There are several important stochastic factors which will be considered in this study.

Assume that the current stage is $m$. $\Delta_m$ is used to denote the information factors to be realized during the time interval $(t_m, t_{m+1}]$ where $t_m$ and $t_{m+1}$ are the time instants corresponding to stages $m$ and $m+1$, respectively. $\Delta_m$ is an abstract conception used to indicate the required information for the state transition of a bus line. Generally speaking, this information always has some random characteristics. When the bus line system evolves into the stage $m+1$ from $m$, the realized information indicated by $\omega_m$ is called a sample of $\Delta_m$.

The uncertainty of $\Delta_m$ mainly comes from two aspects. The first aspect is the passengers' generation at all the stops during $(t_m, t_{m+1}]$. For any given stop, the passengers come to this stop during $(t_m, t_{m+1}]$ with different arrival times, different destination stops and specified behavior types. Many stochastic factors will influence the result of the above passenger generation process. The second aspect is about the uncertain influences of other road traffic and the signal control at intersections on the cruising of buses during $(t_m, t_{m+1}]$. The above influence will reflect in the variable travel times of road segments and the uncertain delays of buses at the signalized intersections.

### 3.3. State Variable

It is hard to follow and supervise a complex system continuously. So a practical way is to specify some time points and only to supervise the system at these time points. In this study, we choose all the starting times of holding operations as the monitoring time points. The total number of stages included in a typical observation period is countable.

The state variable of a bus line system is the chosen features of the system corresponding to the stages. Due to the complicity of a bus line system, it is very difficult to extract a group of features to describe the whole system completely. If the chosen features are too detailed and abundant, the incurred computation will be unbearable later. So we choose only three features as the components of the state variable.

Assume that the current stage is $m$. The first feature chosen to construct the state variable is the latest arrival times of buses at all the stops. We denote this feature with respect to stop $e$ and stage $m$ by $t_e^{@,m}$. The second feature is the time remaining for a bus to be activated the next time. Note that "to be activated" means to consider if a holding operation with specified holding time should be applied to a bus at a stop when both the necessary alighting process for the onboard passengers and the boarding process for the passengers who have arrived at the stop before the bus arrives end. We denote the second feature with respect to bus $b$ and stage $m$ by $t_b^{w,m}$. The third feature is the bus stops where buses will be activated the next time. We denote the third feature with respect to bus $b$ and stage $m$ by $e_b^{w,m}$. Organize these three features in a column

vector $s_m$ as follows

$$s_m \equiv (t_1^{\mathscr{A},m}, t_2^{\mathscr{A},m}, \cdots, t_{n_E}^{\mathscr{A},m}, t_1^{\mathscr{W},m}, t_2^{\mathscr{W},m}, \cdots, t_{n_B}^{\mathscr{W},m}, e_1^{\mathscr{W},m}, e_2^{\mathscr{W},m}, \cdots e_{n_B}^{\mathscr{W},m})^{\perp}. \quad (5)$$

$s_m$ is the state variable corresponding to stage $m$. The superscript "$\perp$" stands for the transposition operation. $E$ is the set of all the bus stops. $n_E$ is the number of bus stops. All the available values of state variable constitute the state space $S$. Obviously, the dimension of $S$ is $n_E + 2n_B$.

One remarkable feature of the above state variable is its completeness when being used to describe the evolution of bus line system with the help of the actually chosen actions. In the following, we present a brief explanation about this feature. To plot the trajectory of a given bus, we need to know three key time points with respect to a bus stop consisting of the arrival time of this bus, the starting time of holding operation and the departure time of this bus. In the state variable, the arrival times have been noted down. If the current system time is $t$, then we know the starting time of holding operation is $t + t_b^{\mathscr{W},m}$. Suppose that the bus $b$ will be activated with holding time $a(m)$ at stage $m$. The departure time of bus $b$ can be calculated by $t + t_b^{\mathscr{W},m} + a(m)$. It is easy to plot the whole trajectory of any given bus with the above three time attributes.

### 3.4. Policy

A policy is a mapping from the state space to the action set. For a given value of the state variable, there is a unique action according to a policy. Denote a policy by $\pi \in \Pi$ where $\Pi$ is the set of all available policies. The above relationship between states and actions can be expressed by

$$a = A^{\pi}(s). \quad (6)$$

The aim of ADP is to determine the optimal policy satisfying some given constraints. In the context of our application, a policy is corresponding to a strategy for determining the optimal holding time of a bus at a stop. The idea underlying the above policy optimization will be clear after the introduction of state value function in the next subsection.

### 3.5. State Value and Q-factor

When an action is determined, its impact on the stability of bus line system needs to be measured so as to judge the quality of the action. Inspired by some existing studies (Daganzo 2009b; Bartholdi and Eisenstein 2012; He 2015), we propose two kinds of headways to be used to construct the proper measures of the above impact. One is the Expected System Headway (**ESH**); the other is the Dynamic Circle Headway (**DCH**).

The ESH is defined by

$$\tilde{H} = (X + \sum_{i \in I} t_i^W \tilde{v} + \sum_{e \in E} t_e^D \tilde{v}) / (n_B \tilde{v}). \quad (7)$$

In eq. (7), $X$ is the length of the bus line; $\tilde{v}$ is the average bus cruising speed; $t_i^W$ is the expected delay for a bus at intersection $i \in I$ where $I$ is the set of all the intersections; $t_e^D$ is the expected dwell time for a bus at stop $e$. The DCH has been defined by eq. (1) in subsection 2.2.

With these two headways, we can define the cost of action $a$ as follows

$$c_a = \sum_{b \in B}(h_b - \mathcal{H})^2, \tag{8}$$

where the headway $h_b$ is corresponding to the time when control action $a$ ends. $\mathcal{H}$ stands for ESH or DCH. The cost $c_a$ of action $a$ with DCH as the coefficient $\mathcal{H}$ has the same changing trend as $\sigma_H(t)$ defined in eq. (2).

The related computation to obtain ESH is carried out once for all. But DCH needs to be calculated at every stage. Both ESH and DCH can be used in the Q-learning algorithm to be given. The numerical results about to be presented in section 5 show that DCH is more effective than ESH in the same situation. But ESH has its own use in many other methods, e.g. the terminal station holding control, which will be considered in the numerical experiments later.

For any stage of the bus line system, there is a corresponding action. If we knew the total number of stages during an observation period in advance, we could add up all $c_a$ with respect to these stages. To minimize the sum obtained above could be used to find out the optimal actions for all the stages, in another word to determine the optimal holding strategy. Generally speaking, it is impossible to know the exact number of stages to be encountered during the observation period in advance, let alone to calculate the exact value of $c_a$ in advance. By the way if the observation period is lengthened, the total number of stages will be considerably large. In view of the above analysis, it is better to use a discounted cost of action $a$ to replace the undiscounted cost $c_a$ in the objective. The modified objective function is as follows

$$\min_{\pi \in \Pi} \mathbb{E}\{\sum_{m=0}^{\infty}[\gamma^m c(s_m, A^\pi(s_m), \Delta_m)]\}, \tag{9}$$

where $\gamma < 1$ is the discounted rate and $\Delta_m$ is the information factors to be realized during the time interval from stage $m$ to $m+1$. Note that $c(s_m, A^\pi(s_m), \Delta_m)$ is the cost of action $A^\pi(s_m)$ with respect to stage $m$. This cost can be calculated during the observation period using Eq. (8). To guarantee that the value of the objective function is finite, we assume that the cost of any feasible action is smaller than a given positive value. The operator $\mathbb{E}$ is to obtain the expected value of some random variable.

Now we can define the optimal policy by

$$\pi^* = \arg\min_{\pi \in \Pi} \mathbb{E}\{\sum_{m=0}^{\infty}[\gamma^m c(s_m, A^\pi(s_m), \Delta_m)]\}. \quad (10)$$

The various constraints is implied in the feasible policy set $\Pi$. It is nearly impossible to solve the above optimization problem directly. Fortunately, the Bellman Equation saves us by supplying an equivalent optimal condition with respect to any two consecutive stages. The cost-to-go or the value of being in a state $s$ is denoted by $V(s)$. Use the superscript "*" to indicate the optimal value. According to the Bellman equation, the optimal $V(s)$ should satisfy the following equation

$$V^*(s_m) = \mathbb{E}[c(s_m, A^{\pi^*}(s_m), \Delta_m) + \gamma V^*(s_{m+1} | s_m, A^{\pi^*}(s_m), \Delta_m)]. \quad (11)$$

Before to present the modified Q-learning algorithm, we should know the basic relationship between state value $V(s)$ and Q-factor $Q(s,a)$ that is as follows

$$V(s) = \min_a Q(s,a). \quad (12)$$

Comparing with $V(s)$, the dimension of the independent variables in $Q(s,a)$ is enlarged by one. The adoption of $Q(s,a)$ brings us the convenience to choose the optimal action at state $s$ by directly comparing the values of $Q(s,a)$ over action set $A$. On the contrary, if $V(s)$ is used, in order to obtain the optimal action we need to push the system forward with a chosen action and consider all the following states. The above process will be very complicated.

The abstract form of Q-factor $Q(s,a)$ must be replaced by a concrete one; or else, the algorithm about to be presented cannot be implemented. To overcome the above difficulty, we can substitute this abstract function for some approximate function $Q(s,a,\vartheta)$ with undetermined parameter vector $\vartheta$. A popular tool commonly used to deal with this situation in the ADP realm is the Artificial Neural Network (ANN).

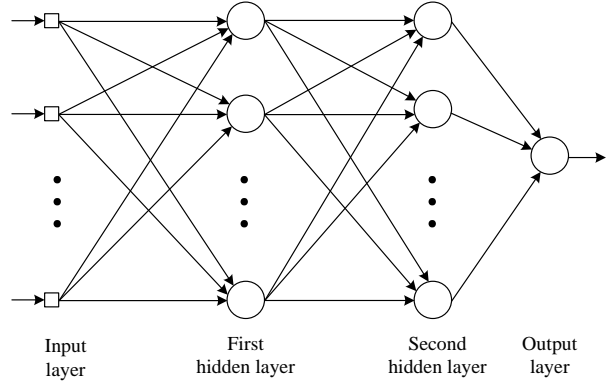

Figure 1 Multilayer perceptron with two hidden layers and one output neuron.

In this study we will use a multilayer perceptron with two hidden layers and one output neuron as shown in Figure 1 to approximate the Q-factor. The number of nodes in the input layer is $n_E + 2n_B + 1$ because the input includes the elements of state variable $s$ and the action $a$. Every connector (or synapse in terms of ANN) has an undetermined weight. Every node (or neuron in terms of ANN) in the hidden and output layers has an undetermined bias. These undetermined weights and biases constitute the undetermined parameter vector $\vartheta$. Every node in the hidden and output layers also has its sigmoid function. The logistic function is defined by $\varphi(v) = [1+\exp(-\alpha v)]^{-1}$ where $\alpha$ is the slope parameter of the sigmoid function. We will choose the above logistic function as the sigmoid function actually used in our study.

For every node in the input layer, its output signal is equal to its input signal that is an element of state variable or the action $a$. For every node in the hidden and output layers, the process of obtaining its output signal includes the following four steps. Firstly, for every upstream node connected to the current node, we multiply the output signal of this upstream node by the weight of the related connector. Secondly, we sum up the products obtained above. Thirdly, we subtract the bias of the current node from the sum obtained above. At last, we use the difference obtained above as the input variable $v$ to compute the value of the sigmoid function $\varphi(v)$. This value is the output signal of the current node. The output signal of the only node in the output layer is the estimated value of Q-factor. Since we only consider the standard operations of ANN in this study, readers can refer to textbooks about ANN, such as Haykin (2009), for a more detailed explanation.

## 4. Q-learning Algorithm with Multistage Look-ahead

To make a concise presentation of the modified Q-learning algorithm, we will first introduce the process of the state transition and then deal with choosing action by the multistage look-ahead mechanism. At last the concrete process of the modified Q-learning algorithm will be presented.

### 4.1. State Transition Process

Assume the current stage is $m$. In order to generate the new state $s_{m+1}$ given state $s_m$, action $a_m$, information sample $\varphi_m$ and the current system time $\vec{t}$, the following process should be carried out.

**Step 1**: Find out the bus to be activated. Here a bus $b \in \arg\min_{i \in B}\{t_i^{\mathbb{U},m}\}$ will be chosen to be activated. If there is more than one bus satisfying the condition above, the choice can be arbitrary. Let $t_{\min} := \min_{i \in B}\{t_i^{\mathbb{U},m}\}$.

**Step 2**: Activate the chosen bus $b$ using the information sample $\omega_m$ of $\Delta_m$. Here we assume that bus $b$ will depart from its current stop $e$ and enter the bus line segment $g$. Renew the attributes $t_b^{\mathbb{U},m}$ with $t_b^{\mathbb{U},m} := t_b^{\mathbb{U},m} + a_m + t_g^{b,m} + t_{e\oplus 1}^{b,m}$ and $e_b^{\mathbb{U},m}$ with $e_b^{\mathbb{U},m} := e_b^{\mathbb{U},m} \oplus 1$, respectively. Note that $e \oplus 1$ is the nearest downstream bus stop of stop $e$. $t_g^{b,m}$ is the travel time of the bus line segment $g$ for bus $b$ following the stage $m$. $t_{e\oplus 1}^{b,m}$ is the forced dwell time of bus $b$ at stop $e \oplus 1$. $t_{e\oplus 1}^{b,m}$ is equal to the bigger one between the alighting time for passengers who get off bus $b$ at stop $e \oplus 1$ following the stage $m$ and the boarding time required by passengers who have arrived at $e \oplus 1$ before bus $b$ arrives at stop $e \oplus 1$.

**Step 3**: Renew the related components of state variable as follows

$t_{e\oplus 1}^{\mathcal{A},m+1} := t_b^{\mathcal{W},m} - t_{e\oplus 1}^{b,m} + \vec{t}$ and $t_\theta^{\mathcal{A},m+1} := t_\theta^{\mathcal{A},m}, \forall \theta \neq e$;

$t_i^{\mathcal{W},m+1} := t_i^{\mathcal{W},m} - t_{\min}$ and $e_i^{\mathcal{W},m+1} := e_i^{\mathcal{W},m}, \forall i \in B$.

Now with the renewed components of $t_\theta^{\mathcal{A},m+1}$, $t_i^{\mathcal{W},m+1}$ and $e_i^{\mathcal{W},m+1}$, the next state $s_{m+1}$ can be constructed.

The above state transition process is used in the off-line version of the Q-learning algorithm to mimic the operation of a real bus line. The information sample $\varphi_m$ is used to produce the sample values of $t_g^{b,m}$ and $t_{e\oplus 1}^{b,m}$. If the on-line Q-learning is carried out, the sample value of $t_g^{b,m}$ will be replaced by its real value and the sample value of $t_{e\oplus 1}^{b,m}$ will be deduced from the actually arrived passengers at stop $e\oplus 1$ and the passengers riding in bus $b$. To simplify the subsequent expressions, we will use an abstract function $s_{m+1} = S^\Theta(s_m, a_m, \varphi_m)$ to represent the above state transition process.

**4.2. Choosing an Action by Multistage Look-ahead**

To choose an action by looking ahead several stages is to roll forward several stages of the simulation system to see which action will bring us the best result. To roll the simulation system forward, we need to use the expected values for many quantities, such as the travel times of bus line segments and the dwell times of buses. These expected values are different from the actual or sample values that we have used in the state transition process in the preceding subsection.

Assume that the current time is $\vec{t}$, the current system state is $s_m$ and the number of successive stages to be looked ahead is $N$. The explanations of some new notations are as follows. $e_b^{w,i}$ stands for the bus stop for bus $b$ to be activated when the current looking ahead stage is $i \in \{1,2,\ldots N\}$. $t_b^{w,i}$ is the time remaining for bus $b$ to be activated at the looking ahead stage $i \in \{1,2,\ldots N\}$. $t_e^{la,i}$ is the modified latest arrival time corresponding to bus stop $e$ when the current looking ahead stage is $i \in \{1,2,\ldots N\}$. Note that we use the time point $\vec{t}$ as the origin to define the modified latest arrival time. To distinguish the expected value from the actual or sample value, the over bar will be used to indicate the expected values. So $\bar{t}_g^b$ stands for the expected travel time of bus line segment $g$ for bus $b$. $\bar{t}_e^b$ is the expected dwell time of bus $b$ at stop $e$.

The process of searching the optimal action $a^*$ with $N$-stage is as follows.

**Step 1**: Assign a large positive value to $\tilde{c}_i, \forall i \in \{1, 2, \cdots, N\}$.

**Step 2**: Choose a bus $b^1 \in \arg\min_{d \in B}\{t_d^{\tilde{0},m}\}$ to be active at the first level.

**Step 3:** Carry out the following multiple nested "for" loops.

    **For** every action $a^1 \in A_{e_{b^1}^{\tilde{0},m}}$, execute the following operations (start the first level):

$$e_{b^1}^{w,1} := e_{b^1}^{\tilde{0},m} \oplus 1 \text{ and } e_b^{w,1} := e_b^{\tilde{0},m}, \forall b \neq b^1;$$

$$t_{b^1}^{w,1} := t_{b^1}^{\tilde{0},m} + a^1 + \overline{t}_{g^1}^{b^1} + \overline{t}_{e_{b^1}^{w,1}}^{b^1} \text{ and } t_b^{w,1} := t_b^{\tilde{0},m}, \forall b \neq b^1;$$

$$t_{e_{b^1}^{w,1}}^{la,1} := t_{b^1}^{\tilde{0},m} + a^1 + \overline{t}_{g^1}^{b^1} \text{ and } t_e^{la,1} := t_e^{\partial,m} - \vec{t}, \forall e \neq e_{b^1}^{w,1}.$$

Note that the two ends of bus line segment $g^1$ are $e_{b^1}^{\tilde{0},m}$ and $e_{b^1}^{w,1}$.

Calculate the cost of action $c_{a^1}$ by eq. (8).

Choose a bus $b^2 \in \arg\min_{d \in B}\{t_d^{w,1}\}$ to be activated at the second level.

**For** every $a^2 \in A_{e_{b^2}^{w,1}}$, execute the following operations (start the second level):

$$e_{b^2}^{w,2} := e_{b^2}^{w,1} \oplus 1 \text{ and } e_b^{w,2} := e_b^{w,1}, \forall b \neq b^2;$$

$$t_{b^2}^{w,2} := t_{b^2}^{w,1} + a^2 + \overline{t}_{g^2} + \overline{t}_{e_{b^2}^{w,2}}^{b^2} \text{ and } t_b^{w,2} := t_b^{w,1}, \forall b \neq b^2;$$

$$t_{e_{b^2}^{w,2}}^{la,2} := t_{b^2}^{w,1} + a^2 + \overline{t}_{g^2} \text{ and } t_e^{la,2} := t_e^{la,1}, \forall e \neq e_{b^2}^{w,2}.$$

Calculate the cost of action $c_{a^2}$ by eq. (8).

…

Choose a bus $b^N \in \arg\min_{d \in B}\{t_d^{w,N-1}\}$ to be activated at the level $N$.

**For** every $a^N \in A_{e_{b^N}^{w,N-1}}$, execute the following operations (start the level $N$):

$$e_{b^N}^{w,N} := e_{b^N}^{w,N-1} \oplus 1 \text{ and } e_b^{w,N} := e_b^{w,N-1}, \forall b \neq b^N;$$

$$t_{b^N}^{w,N} := t_{b^N}^{w,N-1} + a^N + \overline{t}_{g^N} + \overline{t}_{e_{b^N}^{w,N}}^{b^N} \text{ and } t_b^{w,N} := t_b^{w,N-1}, \forall b \neq b^N;$$

$$t_{e_{b^N}^{w,N}}^{la,N} := t_{b^N}^{w,N-1} + a^N + \overline{t}_{g^N} \text{ and } t_e^{la,N} := t_e^{la,N-1}, \forall e \neq e_{b^N}^{w,N}.$$

Calculate the cost of action $c_{a^N}$ by eq. (8);

Add the estimate of the future states' value to $c_{a^N}$ as follows

$$\vec{c}_{a^N} := c_{a^N} + \gamma \min_{a_{N+1}} Q(s_{N+1}, a_{N+1}, \vartheta).$$

Note that $s_{N+1}$ and $a_{N+1}$ are the state and an available action after $a^N$ has been carried out, respectively. $\vec{c}_{a^N}$ is the sum of $c_{a^N}$ and the estimate of the future states' discounted value.

Let $\tilde{c}_N := \min\{\tilde{c}_N, \vec{c}_{a^N}\}$.

**End** the "for" of the level $N$.

Let $\tilde{c}_{N-1} := \min\{\tilde{c}_{N-1}, c_{a^{N-1}} + \gamma \tilde{c}_N\}$.

…

Let $\tilde{c}_2 := \min\{\tilde{c}_2, c_{a^2} + \gamma \tilde{c}_3\}$.

**End** the "for" of the level 2.

Let $\tilde{c}_1 := \min\{\tilde{c}_1, c_{a^1} + \gamma \tilde{c}_2\}$ and $a^* := \arg\min_{a^1}\{\tilde{c}_1, c_{a^1} + \gamma \tilde{c}_2\}$

**End** the "for" of the level 1.

**Step 4**: Output the optimal action $a^*$ with respect to the current state $s_m$.

### 4.3. Process of the Modified Q-learning Algorithm

Generally the classical Q-learning algorithm can pick out the best action easily by using the Q-factor to estimate the value of any available action. But in our case, the Q-factor is approximated by ANN. In the early iterations of Q-learning algorithm, the estimate from ANN with many undetermined parameters is so crude and inaccurate that it is a little help using the estimate to choose an effective action to stabilize the bus line. When the bus line system is large-scale and the number of undetermined parameters is big, it will be very difficult to obtain a good estimate from the early iterations of the original Q-learning algorithm.

To improve the accuracy of the estimate of Q-factor, we introduce the multistage look-ahead mechanism into the classical Q-learning algorithm. The underlying thought is that through looking ahead several successive stages when a proper holding time needs to be determined, we can accelerate the information accumulation with the approximate Q-factor in an ANN form. Especially at the early phase of the implementation of the algorithm, a relatively big number of successive stages to be looked ahead will increase the chance of choosing an effective action and at the same time renewing the undetermined coefficients of ANN along the right directions. The introduction of the multistage look-ahead mechanism can make the algorithm get through the mire of the early iterations more quickly and easily.

Assume that the observation period spans a time period of $[0, T]$. In the following algorithm, $m$ is used to indicate the sequence number of the current stage during an observation period. $K$ is the total number of observation periods to be studied in the algorithm. $k$ is the sequence number of the current observation period in the total $K$ periods. $n_M^k$ is the total number of stages in the observation period $k$. $A_m^k$ is the set of all available actions at stage $m$ in the observation period $k$. $\vartheta$ is the vector of undetermined coefficients, i.e. the weights of connecters and the biases of nodes in ANN. Assume that the number of successive stages to be looked ahead is $\Gamma$.

The detailed process of the modified Q-learning algorithm is as follows.

**Step 1**: Initialize the coefficients $\vartheta$ of the ANN which is used to approximate the Q-factor $Q(s,a)$. The approximated function by ANN is denoted by $\bar{Q}(s,a,\vartheta)$. Set $k:=1$.

**Step 2:** Set the current time $\vec{t}:=0$ and let $\varepsilon:=\varepsilon-k\xi$. Set $m:=0$. Initialize the starting state variable $s_m^k$ by using the information of bus line system at time 0. If a simulation system is used to stand for a bus line, the passengers can be generated at all the bus stops in this step in advance.

**Step 3:** While the condition $\vec{t}\leq T$ holds, do the following "while" loop:

Step 3a. Determine the action using $\varepsilon$-greedy rule. With probability $\varepsilon<1$, choose an action $a_m^k$ at random from $A_m^k$. With probability $1-\varepsilon$, choose $a_m^k$ using $\Gamma$ multistage look-ahead technique.

Step 3b. Compute the next state $s_{m+1}^k = S^\Theta(s_m^k, a_m^k, \omega_m^k)$.

Step 3c. Compute $\hat{q}_m^k = c(s_m^k, a_m^k) + \gamma \min_{a_{m+1}\in A_{m+1}^k} \bar{Q}_{m+1}^{k-1}(s_{m+1}^k, a_{m+1}, \vartheta)$. Here $c(s_m^k, a_m^k)$ will be calculated by function (8) with the given state $s_m^k$ and action $a_m^k$.

Step 3d. Update the coefficients $\vartheta$ to refine the approximate $\bar{Q}_m^k$.

Update $\vartheta$ as follows

$$\vartheta := \vartheta + \lambda \nabla \bar{Q}_m^{k-1}(s_m^k, a_m^k, \vartheta)(\hat{q}_m^k - \bar{Q}_m^{k-1}(s_m^k, a_m^k, \vartheta)).$$

Or we can update the $\vartheta$ using

$$\vartheta := \vartheta + \lambda \nabla \bar{Q}_m^{k-1}(s_m^k, a_m^k, \vartheta)(c(s_m^k, a_m^k) + \gamma \bar{Q}_{m+1}^{k-1}(s_{m+1}^k, a_{m+1}, \vartheta) - \bar{Q}_m^{k-1}(s_m^k, a_m^k, \vartheta)).$$

In the above operation, action $a_{m+1}$ can be decided by the $\varepsilon$-greedy rule.

With the renewed $\vartheta$, the Q-factors can be estimated as follows

$$Q_m^k(s_m^k, a_m^k) = (1-\alpha_{k-1})\bar{Q}_m^{k-1}(s_m^k, a_m^k, \theta) + \alpha_{k-1}\hat{q}_m^k.$$

Step 3e. With the time interval $t_{m,m+1}$ between stage $m$ and $m+1$, set $\vec{t}:=\vec{t}+t_{m,m+1}$. Update $m:=m+1$.

**Step 4:** Let $k:=k+1$. If $k\leq K$, go to step 2.

**Step 5:** Return the Q-factor estimates $(Q_m^K)_{m=1}^{n_M^K-1}$ and $\vartheta$.

In the above algorithm, we adopt the Epsilon-greedy rule to choose an appropriate action, i.e. to determine the length of holding time. The Epsilon-greedy rule is a very popular tool used in the ADP realm. The aim of using this tool in Q-learning algorithm is to exploit the accumulated knowledge about the state value and at the same time to explore the new states that have not been visited. The value of parameter $\varepsilon$ is limited to the range of [0, 1). $\xi$ is a small positive real number used to gradually diminish the value of $\varepsilon$ such that $\varepsilon - K\xi \geq 0$. At the early phase of the implementation of the modified Q-learning algorithm, a relatively big value of $\varepsilon$ can be used; later it will be reduced to a relatively small one. To do so is to put more stress on exploration at the beginning. With the gradually increased $k$, we can make use of the already refined approximate function of Q-factor with ANN form to choose an effective action by decreasing the probability of choosing actions randomly.

$\bar{Q}(s,a,\vartheta)$ is the approximate function explained at the end of the subsection 3.5. With the given state variable $s$ and control action $a$ as the input, the ANN will generate the output $\bar{Q}(s,a,\vartheta)$ with $\vartheta$ as its parameters. If the values of $s$ and $a$ are fixed, we can view $\bar{Q}(s,a,\vartheta)$ as a function of variable $\vartheta$. In the above algorithm, we use $\nabla \bar{Q}(s,a,\vartheta)$ to stand for the partial gradient of $\bar{Q}(s,a,\vartheta)$ with respect to $\vartheta$. Readers can refer to textbooks, such as Haykin (2009), to know the way of computing $\nabla \bar{Q}(s,a,\vartheta)$.

## 5. Numerical Experiments
### 5.1. Simulation System of a Bus Line
A detailed simulation system will be used to verify the effectiveness of our Q-learning algorithm. The key elements of this simulation system will be introduced below.
### 5.1.1. Passenger Generation
Assume that there are two types of passengers characterized by their specified average boarding and alighting times. One type is denoted by $P^s$ with the average boarding time $t^s_{Bd}$ and the average alighting time $t^s_{Al}$; the other type is denoted by $P^q$ with the average boarding time $t^q_{Bd}$ and the average alighting time $t^q_{Al}$. The rate of the number of passengers in $P^s$ over the number of passengers in $P^q$ is $r_{s/q}$.

The passenger generating rate at bus stop $e$ is denoted by $r_e$. During an observation period, passengers will be generated at all bus stops one time unit after the other by the Monte-Carlo method. For any bus stop, there is a set of bus stops which are the possible destination stops for the passengers generated at this stop. Based on field data, a specified probability will be assigned to a specified stop in the above set of stops. Passengers will choose their destination stops according to the given probability distribution.
### 5.1.2. Travel Time

The part of a bus line between two neighboring stops is defined as bus line segment. A bus line segment generally consists of several road segments and several intersections. A road segment is the part of bus line segment between two neighboring ends. Here an end indicates an intersection or a stop. A typical road segment is denoted by $d$. The length of road segment $d$ is denoted by $l_d$. With the average cruising speed $\bar{v}$, the expected travel time $\bar{t}_d$ of $d$ is $l_d/\bar{v}$. The sample travel time of $d$ can be generated by $t_d = \bar{t}_d + \delta_d$ where $\delta_d$ is a normal random variable with zero mean and variance $\sigma_d^2$.

Assume that the signal control at intersection $i$ is the pre-timed two phases with respect to the approach of the bus line in question. The traffic signal cycle length $t_i^{cl}$, the red phase $t_i^{red}$, the green phase $t_i^{green}$, and the remaining time of the initial phase $t_i^{or}$ are given. Here we assume that the amber phase has been appropriately included in the other phases. So the amber phase will be omitted in the following analysis and $t_i^{cl} = t_i^{red} + t_i^{green}$, $\forall i \in I$ holds.

When a bus arrives at an intersection and the current traffic light is green, the actual delay at this intersection is zero; if the current traffic light is red at this time, the actual delay at this intersection equals the remaining time of the red phase. The expected delay at intersection $i$ is given by $\frac{1}{2}(t_i^{red})^2 / t_i^{cl}$.

A typical bus line segment is denoted by $g \in G$ where $G$ is the set of all bus line segments. The expected travel time of $g$ is the sum of all the expected travel times of its road segments and all the expected delays of its intersections. The actual travel time of $g$ with entering time $t_g^{enter}$ is equal to the sum of the actual travel times of related road segments and actual delays of related intersections for a bus which enters the first road segment of $g$ at time $t_g^{enter}$.

**5.1.3. Dwell Time**

Since we assume that boarding and alighting take place at different doors for any bus, the dwell time should be the bigger one between the boarding time and the alighting time.

The sample boarding time is determined as follows. First, when the bus arrives at the stop, the passengers waiting outside start to board the bus. After the above boarding, the holding operation is carried out to the bus. During the holding period, all the passengers who arrive at the current stop before the end of holding are allowed to board the bus. The total boarding time of the above two stages of boarding will be the sample boarding time. Note that the boarding process should be prohibited if the maximal capacity of the bus is reached. The value is obtained by a detailed simulation process. Similarly, the sample alighting time is obtained using the simulation.

When a control action needs to be determined by the multistage look-ahead, the expected dwell time is required. The expected dwell time will be the bigger one between the expected boarding time and the expected alighting time. The expected boarding time will be determined by the time interval between two latest arrival times with respect to the stop, the passenger generation rate and passengers types. The expected alighting time needs to be estimated at the stops which locate between the position of the current stop and the position of the leading bus. When the current bus arrives at one of the above stops, some passengers will alight. Among these alighting passengers, some of them have boarded the bus before the bus departs from the current stop. Some of them are just added after the bus departs from the current stop. The number of the added passengers will be estimated from the passenger arrival rate and the related time interval.

### 5.1.4. Assign Values to Basic Coefficients

The commonly used coefficients in this section are given as follows. The length of a simulation period $T$ is set to 7200 seconds (s). The average cruising speed $\tilde{v}$ is set to 30km/hr. The rate $r_{s/q}$ is set to 1/9. The average boarding times $t_{Bd}^{s}$ and $t_{Bd}^{q}$ are set to 4s and 1s, respectively. The average alighting times $t_{Al}^{s}$ and $t_{Al}^{q}$ are set to 2s and 0.5s, respectively. The discount rate $\gamma$ is set to 0.5. The total number of simulations $K$ is set to 300. If the length of a road segment $d$ is $l_d$ (km), we assume that the related $\delta_d$ is a normal random variable with zero mean and standard deviation $\sigma_d = 5l_d$ (s).

The following coefficients will be used commonly unless otherwise stated. The initial Epsilon $\varepsilon$ is equal to 0.6 and the related reduction $\xi$ is set to 1/600. The unit holding time $\tau$ is set to 2s. The corresponding action set $A$ used at all stops is {0, 2, 4, 6, 8, 10}. To simplify the following expressions, action set is denoted by $A_{n \times m}$. Here the subscript $n \times m$ means that $\tau = n$ (s) and the value of the biggest action equals the product of $n$ times $m$. So the action set {0, 2, 4, 6, 8, 10} can be denoted by $A_{2 \times 5}$. The Dynamic Circle Headway (DCH) is commonly used in this section unless otherwise stated.

The numbers of nodes in the first and second hidden layers of ANN are 5 and 3, respectively. The slope parameter of the sigmoid function is set to 0.5. The initial value of any element of $\vartheta$ is generated randomly from a uniform distribution over [-2, 2]. Unless otherwise stated, the unit of time used in the following tables and figures is always second (s) and the unit of length is always meter (m).

### 5.2. Basic Data of Bus Lines

In this section, we will consider 5 circular bus lines. The basic data about these 5 bus lines are presented in Table 1. These bus lines are different from each other in aspects including the number of buses $n_B$, the number of bus stops $n_E$, the number of intersections $n_I$ and the length $X$.

Table 1

The basic data about 5 bus lines.

| Bus Line | L1 | L2 | L3 | L4 | L5 |
|---|---|---|---|---|---|
| $n_B$ | 5 | 7 | 9 | 11 | 13 |
| $n_E$ | 18 | 24 | 30 | 36 | 42 |
| $n_I$ | 8 | 11 | 13 | 15 | 18 |
| $X$ (km) | 10.65 | 14.21 | 17.95 | 21.35 | 24.6 |

Every component of a bus line has its own features. To give the complete information about all the components of the above 5 bus lines will take up unbearable space of this paper. In view of this situation, we choose only to present the complete information about bus line L5 below. This bus line will be used as the main bus line in the following analysis. Readers can use the information to be presented to rebuild this bus line.

Table 2

Bus stops with given arrival rate in bus line L5.

| $r_e$ | Related Stops in L5 |
|---|---|
| 1 | 4,6,8,13,15,16,17,20,21,24,26,28,31,34,36,38,40 |
| 2 | 1,2,5,7,10,12,14,18,23,25,27,32,33,35,37,39,41,42 |
| 3 | 3,9,11,19,29 |
| 4 | 22,30 |

Note. The unit of $r_e$ is passengers per minute in this table.

We assume that there are four types of passenger arrival rates. In Table 2, the sequence numbers of bus stops with given arrival rate in bus line L5 are grouped.

Table 3

The probabilities for the following downstream stops to be chosen as destination.

| No. | 1 | 2 | 3 | 4 | 5 | 6 | 7 | 8 | 9 | 10 | 11 | 12 | 13 |
|---|---|---|---|---|---|---|---|---|---|---|---|---|---|
| Series 1 | 0.0135 | 0.027 | 0.0541 | 0.0811 | 0.1081 | 0.1351 | 0.1351 | 0.1216 | 0.1216 | 0.0811 | 0.0541 | 0.0405 | 0.0270 |
| Series 2 | 0.0345 | 0.0862 | 0.1207 | 0.1552 | 0.1724 | 0.1552 | 0.1207 | 0.0862 | 0.0517 | 0.0172 | / | / | / |

In order to choose a destination stop for a new generated passenger at a stop $e$, we need to know the probabilities of the following downstream stops of stop $e$. Two series of probabilities given in Table 3 will be used in the simulation. Series 1 and 2 have 13 and 10 elements, respectively. The sum of the values of elements in any one series is 1. These two series can be regarded as a discrete probability distribution.

Suppose that the current stop uses series 1 as its series of probabilities to choose destination stops for the passengers generated at this stop. The $n$ th element of series 1 is the probability of the $n$ th downstream stop of the current stop to be chosen as a destination stop. In bus line L5, we assume that series 1 is used by bus stops including 1, 7, 10, 13, 20, 21, 27, 30, 31, 37, 40 and series 2 is used by the other remaining stops.

Table 4

The basic data of intersections in bus line L5.

| No. of Intersection | 1 | 2 | 3 | 4 | 5 | 6 | 7 | 8 | 9 | 10 | 11 | 12 | 13 | 14 | 15 | 16 | 17 | 18 |
|---|---|---|---|---|---|---|---|---|---|---|---|---|---|---|---|---|---|---|
| $t_i^{red}$ (s) | 40 | 40 | 40 | 30 | 30 | 40 | 40 | 30 | 30 | 40 | 40 | 40 | 30 | 40 | 40 | 40 | 30 | 40 |
| $t_i^{green}$ (s) | 50 | 30 | 35 | 45 | 30 | 30 | 45 | 35 | 45 | 50 | 30 | 35 | 45 | 50 | 30 | 35 | 45 | 50 |
| $t_i^{or}$ (s) | 20 | 20 | 10 | 20 | 20 | 20 | 30 | 20 | 20 | 10 | 20 | 10 | 20 | 20 | 15 | 20 | 20 | 20 |
| Initial phase | 2 | 1 | 1 | 2 | 2 | 1 | 2 | 2 | 2 | 2 | 1 | 1 | 2 | 2 | 1 | 1 | 2 | 2 |
| Bus line segment | 1 | 4 | 8 | 10 | 12 | 13 | 16 | 16 | 20 | 21 | 24 | 28 | 30 | 31 | 34 | 38 | 40 | 41 |

The basic data of 18 intersections in bus line L5 is presented in Table 4. The red phase $t_i^{red}$, the green phase $t_i^{green}$, and the remaining time of the initial phase $t_i^{or}$ are listed in three rows with second as the time unit. The row labeled by "Initial phase" specifies the initial phases of all the intersections. In this row, 1 and 2 indicate the red and green phases, respectively. The last row points out the corresponding bus line segment which includes the given intersection.

Table 5

The basic data about the buses cruising on L5.

| No. of Bus | 1 | 2 | 3 | 4 | 5 | 6 | 7 | 8 | 9 | 10 | 11 | 12 | 13 |
|---|---|---|---|---|---|---|---|---|---|---|---|---|---|
| Capacity | 72 | 70 | 80 | 60 | 72 | 60 | 72 | 80 | 60 | 72 | 70 | 70 | 80 |
| Initial stop | 1 | 4 | 8 | 11 | 15 | 18 | 21 | 25 | 28 | 31 | 34 | 37 | 40 |
| TRAB (s) | 20 | 0 | 40 | 30 | 50 | 10 | 30 | 36 | 24 | 18 | 26 | 16 | 34 |

In Table 5, the passenger capacity, the initial stops where buses are dwelling at or will reach first at the beginning of simulation, and the Time Remaining for a bus to be Activated at the Beginning of every simulation (TRAB) are presented. From the data in the row indicated by "Initial stop", we can see that at the beginning of simulation the buses are evenly dispersed along the bus line.

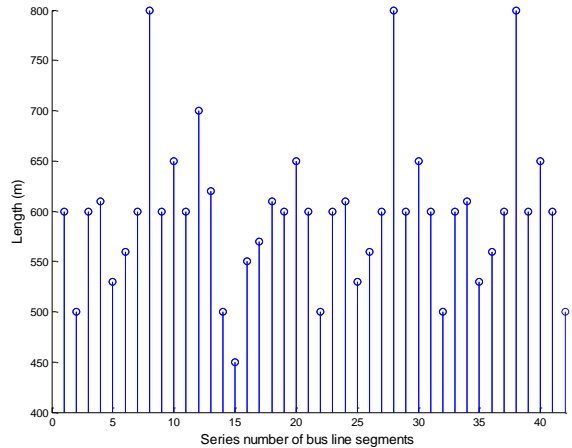

Figure 2 The lengths of bus line segments of L5.

To avoid a very long table, the lengths of bus line segments of L5 are shown in Figure 2. The longest bus line segment is 800 meters (m); the shortest one is 450m. Most of bus line segments are about 600m. The distribution of lengths is relatively random.

**5.3. Main Results Related to the Bus Line with 42 Stops**

In this subsection, we will use bus line L5 as the representative to investigate the performances of various methods. The indices defined in Section 2 will be the main object of investigation here. We will also study the impacts of different action sets and different values of $\varepsilon$ on the simulation results. The performances of two headways used in the cost of action are to be compared at the end of this subsection.

**5.3.1. Introduction of Terminal Station Control**

Before the main result is presented, we need to explain two holding strategies including the single terminal station control and the two terminal stations control. For simplicity, the former is also called the Single Point control (SP) and the latter is also called the Two Points control (TP). The reason why we consider these two holding methods lies in two aspects. On the one hand, these methods are commonly used in practice and have been proved very effective to resist bus bunching sometimes. On the other hand, these methods are very similar to the method proposed by Bartholdi and Eisenstein (2012). The method proposed by Bartholdi and Eisenstein (2012) making use of the backward headway of the current bus at the terminal station has been compared in the existing studies (Xuan et al. 2011; He 2015) with many other methods. The comparing results in these studies show that the performance of this method is relatively good.

Since most existing holding methods adopt some assumptions to simplify the real bus line system, this makes it very hard if not impossible to compare them with the method proposed in this paper in a context of detailed simulation. For example, many existing studies overlooked the impact of intersections with traffic signal control. In view of the second aspect mentioned in the above paragraph, comparing the terminal station control methods with our method supplies an indirect way to compare our method with the other existing methods.

The operation processes of both the single point control and the two points control are the same. The difference only lies in the number of terminal stations used as control point. Their common operation process is as follows. When a bus finishes the boarding and alighting processes at a terminal station used as control point, the headway between this bus and the nearest downstream bus is calculated. Denote this headway by $\vec{H}$. If $\vec{H}$ is longer than the Expected System Headway $\tilde{H}$, the current bus will not be held at this terminal station and departures from this terminal station at once. If $\vec{H}$ is shorter than $\tilde{H}$, the current bus will be held at this terminal station for the time interval of $\tilde{H} - \vec{H}$. The underlying mechanism of these two holding methods is to try to evenly dispatch buses at the specified terminal stations with the expected departure interval of $\tilde{H}$.

**5.3.2. Trajectories' Comparison under Different Control Schemes**

To have an intuitive impression of the operation of the bus line L5, first let us see several figures of the trajectories of buses under different control schemes.

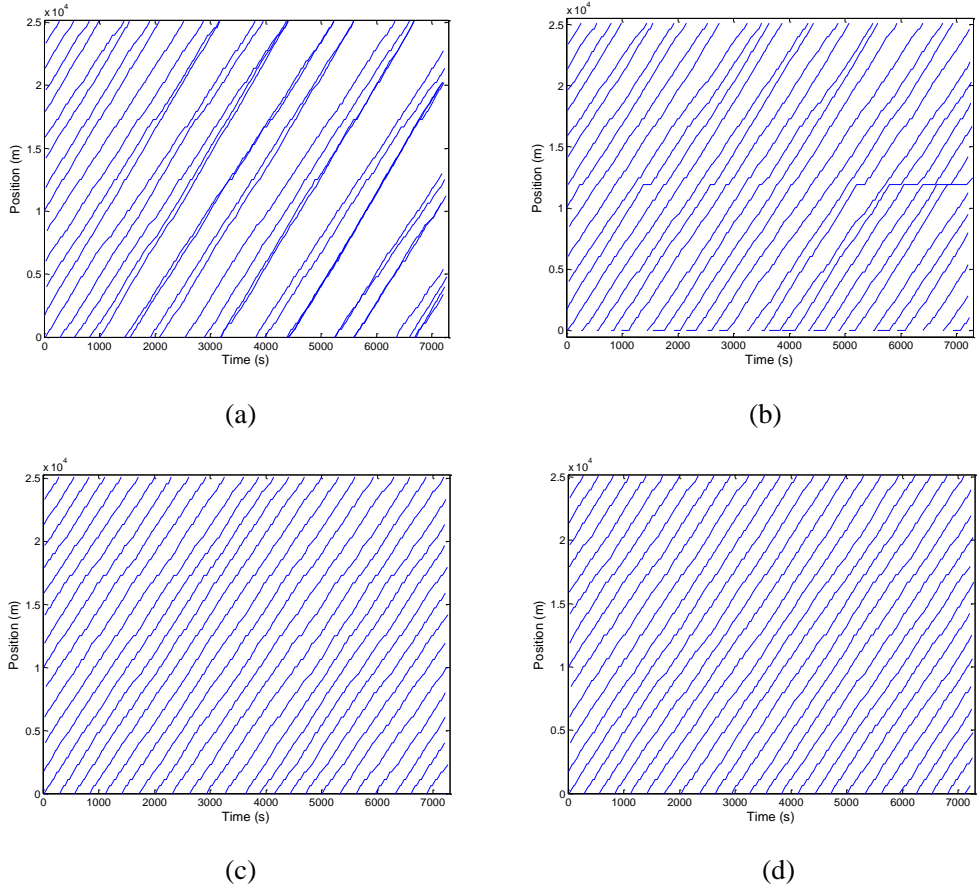

Figure 3. The trajectories of buses in L5: (a) without control; (b) with control at two terminal stations 1 and 21; (c) resulted from QL1S; (d) resulted from QL3S.

The bus line L5 has a strong trend of bus bunching. In a simulation of 2 hours, it is unavoidable for L5 to run into bus bunching. In Figure 3a, a typical simulation result of the trajectories of buses is presented. Though at the beginning of any simulation buses are dispersed quite evenly along the bus line, bunching will appear and become more and more serious with time elapsing.

If we choose stops 1 and 21 as the control points and use the two points control explained above, a typical simulation result of the trajectories of buses is shown in Figure 3b. The expected system headway $\tilde{H}$ is about 275.0 (s). From the trajectories in Figure 3b, we can see that the trajectories of buses have the strong trend to bunch with each other. Fortunately, by holding buses at stops 1 and 21, the whole system runs smoothly and steadily most of the time.

Use "QLnS" to stand for the Q-learning algorithm with the look-ahead of n successive stages. For example, QL4S stands for the Q-learning algorithm with the look-ahead of 4 successive stages. In Figure 3c, the trajectories of buses in L5 resulted from QL1S are plotted. The operation of bus line L5 looks relatively smooth and steady. In Figure 3d, the trajectories of buses in L5 resulted from QL3S demonstrate the similar steady feature as in Figure 3c. Later using the performance indices defined in subsection 2.1, we will distinguish the above intuitively similar trajectories. Note that the results of above two figures come from the last time of 300 simulations and the continuously refined Q-factor during the preceding 299 times simulations is used in the last simulation.

**5.3.3. Learning Process of the Modified Q-learning Algorithm**

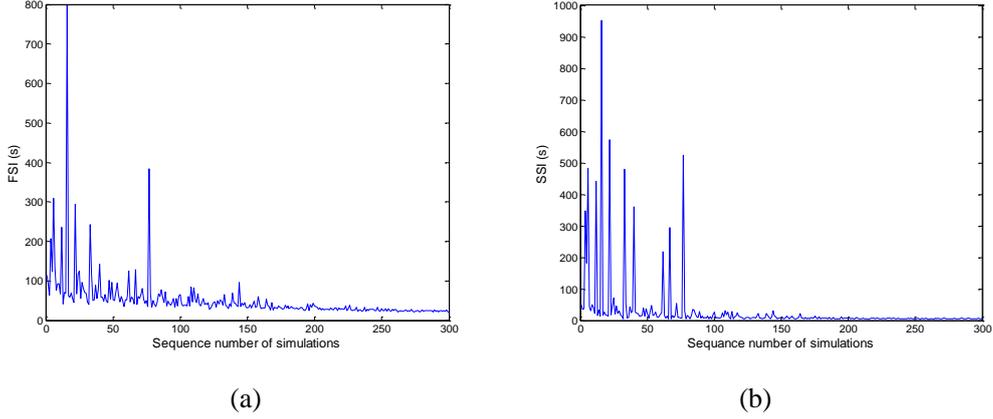

Figure 4 The changing trends with $\varepsilon$ =0.6 and $\xi$ =1/600: (a) the First Stability Index (FSI) $\overline{c}_H$ ; (b) the Second Stability Index (SSI) $\sigma_{\overline{c}}$.

Though in practice we always focus on the final control result of the method that is the performance of bus line system in the last time of simulations in our numerical example, the performance of the method during the whole series of simulations supplies the chance for us to observe the effectiveness and convergence of the Q-learning algorithm.

In Figures 4 the changing trends of the First Stability Index (FSI) $\overline{c}_H$ and the Second Stability Index (SSI) $\sigma_{\overline{c}}$ with $\varepsilon$ =0.6 and $\xi$ =1/600 are plotted, respectively. FSI and SSI fluctuate greatly during the first 100 simulations. With the increasing times of simulations, the range of fluctuation reduces gradually. The big values of $\overline{c}_H$ and $\sigma_{\overline{c}}$ mean the happening of bus bunching in the corresponding simulation. The small values of $\overline{c}_H$ and $\sigma_{\overline{c}}$ mean the smooth and steady operation of bus line system in the corresponding simulation. The changing processes of these two indices demonstrate the learning process of the modified Q-learning algorithm.

Generally, the Epsilon-greedy rule used in the Q-learning algorithm has a great impact on the performance of the algorithm. To demonstrate this impact, we change the values of both $\varepsilon$ and $\xi$ from 0.6 and 1/600 to 0.4 and 1/1000, respectively. The new changing trend of $\overline{c}_H$ in Figure 5a is similar to the one in Figure 4a. The new changing trend of $\sigma_{\overline{c}}$ in Figure 5b is similar to the one in Figure 4b. The main difference between the new results and the old ones lies in the ranges of fluctuations of $\overline{c}_H$ and $\sigma_{\overline{c}}$ especially during the first 100 times of simulations. This difference can be explained from two aspects. On the one hand, a big value of $\varepsilon$ will make the Q-learning algorithm stress on the exploration of various new choices at the beginning and perform better later. On the other hand, a small value of $\varepsilon$ will make the algorithm stress on the exploitation of the known optimal choices at the beginning and take more times of simulations to achieve a better performance later.

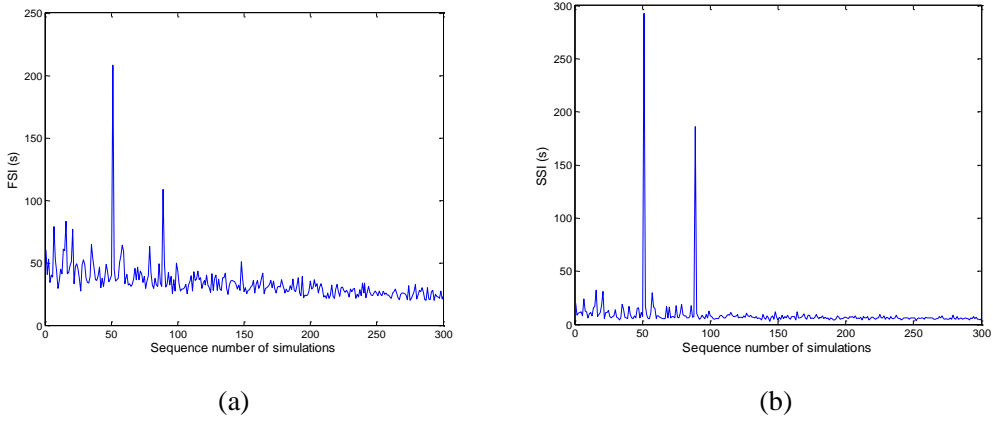

(a)  (b)

Figure 5 The changing trend with $\varepsilon = 0.4$ and $\xi = 1/1000$: (a) the First Stability Index (FSI) $\bar{c}_H$; (b) the Second

Stability Index (SSI) $\sigma_{\bar{c}}$.

**5.3.4. Comparison of Different Control Schemes**

"NC", "SP" and "TP" will be used to indicate the no control, the single terminal station control and the two terminal stations control in the following tables, respectively. "OQL" indicates the original Q-learning algorithm without look-ahead in the following tables. The main results of various methods applied to L5 are presented in Tables 6 and 7.

To verify the effectiveness of the new algorithm, the data to be given in the following tables will be some average values. After running the simulation system 50 times under the no control scheme, the single terminal station control and two terminal stations control, respectively, the average values over 50 times will be used as the results shown in the following tables. The results of various Q-learning algorithms come from the average values over 50 times of simulations. In these 50 times simulations, only the optimal actions are chosen by the multistage look-ahead mechanism. The optimized $\vartheta$ resulted from the Q-learning algorithm with $K = 300$ is used in these 50 times of simulations. Since we allow passengers to board and alight when a bus is held at a stop, part of the holding time will be used to board and alight. The column indicated by "$\mathbb{C}_a$" is the sum of the holding times during which there is no passenger boarding or alighting.

Table 6

The indices of system stability and control action of L5 under different control schemes.

| Methods | $\sum_{t \in \bar{T}} \sigma_H(t)$ | $\bar{c}_H$ | $\sigma_{\bar{c}}$ | $\max_{t \in \bar{T}} \sigma_H(t)$ | $\min_{t \in \bar{T}} \sigma_H(t)$ | $a_\Sigma$ | $\mathbb{C}_a$ | $\bar{a}$ | $\sigma_{\bar{a}}$ | $n_M$ | Bunch |
|---|---|---|---|---|---|---|---|---|---|---|---|
| NC | 322346.23 | 301.54 | 486.98 | 3207.54 | 34.82 | / | / | / | / | / | Yes |
| SP | 107515.30 | 103.08 | 45.82 | 205.97 | 32.73 | 2217 | 2147 | 88.70 | 99.29 | 25 | Yes |
| TP | 72665.97 | 70.96 | 20.60 | 127.38 | 35.03 | 4292 | 4161 | 85.85 | 92.59 | 50 | Yes/No |
| OQL | 29857.09 | 29.18 | 8.97 | 54.55 | 12.47 | 4918 | 4461 | 4.81 | 4.52 | 1023 | Yes/No |
| QL1S | 25074.97 | 24.43 | 5.88 | 44.86 | 11.96 | 4664 | 4214 | 4.57 | 4.49 | 1021 | No |
| QL2S | 24232.67 | 23.73 | 5.42 | 41.31 | 11.91 | 4480 | 4056 | 4.37 | 4.47 | 1026 | No |
| QL3S | 20999.70 | 20.46 | 5.25 | 40.02 | 9.71 | 4366 | 3977 | 4.25 | 4.38 | 1026 | No |
| QL4S | 22002.77 | 21.55 | 5.50 | 41.42 | 10.97 | 4590 | 4173 | 4.50 | 4.50 | 1021 | No |

| QL5S | 23753.48 | 23.31 | 5.71 | 42.21 | 11.53 | 4610 | 4186 | 4.52 | 4.53 | 1019 | No |

The last column of Table 6 with the title "Bunch" is used to indicate whether bus bunching happens during the simulation. "Yes/No" means that bunching may happen sometimes if we run the corresponding algorithm many times. The data show that Q-learning algorithms with multistage look-ahead outperform other methods in view of the indices of system stability and the outside interference intensity. The single point control cannot remove bunching from the operation of the bus line system. The performances of both the two points control and the original Q-learning algorithm are barely satisfactory. With the increasing number of successive stages to be looked ahead, the performance of the corresponding modified Q-learning algorithm is improved at the beginning; but later the trend is reversed. On the one hand, the above observation shows that the modified Q-learning algorithm with multistage look-ahead can stabilize a high frequency bus line effectively. On the other hand, the increasing randomness due to the increasing number of successive stages will undermine the benefit that just comes from the operation of looking ahead more stages to choose a proper control action.

Table 7

The indices of service level of L5 under different control schemes.

| Methods | $P_1$ | | | | | | | $P_2$ | | | | | $P_3$ | | |
|---|---|---|---|---|---|---|---|---|---|---|---|---|---|---|---|
| | $n_{P_1}$ | $t_{P_1}^W$ | $\sigma_{P_1}^W$ | $t_{P_1}^R$ | $\sigma_{P_1}^R$ | $t_{P_1}^{Tr}$ | $\sigma_{P_1}^{Tr}$ | $n_{P_2}$ | $t_{P_2}^W$ | $\sigma_{P_2}^W$ | $t_{P_2}^R$ | $\sigma_{P_2}^R$ | $n_{P_3}$ | $t_{P_3}^W$ | $\sigma_{P_3}^W$ |
| NC | 8067 | 208.5 | 141.5 | 487.8 | 218.6 | 696.3 | 265.3 | 645 | 335.5 | 199.0 | 211.1 | 227.4 | 321 | 179.0 | 145.1 |
| SP | 8361 | 172.7 | 110.3 | 491.1 | 223.5 | 663.8 | 251.1 | 664 | 226.1 | 151.2 | 230.8 | 221.2 | 166 | 116.0 | 95.1 |
| TP | 8415 | 166.3 | 94.2 | 497.4 | 231.7 | 663.7 | 252.2 | 671 | 179.3 | 97.2 | 218.6 | 223.9 | 93 | 64.4 | 47.6 |
| OQL | 8331 | 148.4 | 87.9 | 502.9 | 227.9 | 651.4 | 246.6 | 643 | 151.7 | 83.3 | 227.2 | 214.7 | 87 | 68.1 | 48.4 |
| QL1S | 8318 | 150.0 | 91.9 | 501.7 | 233.3 | 651.7 | 250.0 | 657 | 146.0 | 83.5 | 214.9 | 214.8 | 93 | 81.1 | 51.1 |
| QL2S | 8383 | 148.6 | 89.5 | 502.4 | 225.1 | 651.0 | 242.0 | 647 | 145.3 | 80.9 | 216.4 | 215.2 | 79 | 65.2 | 50.1 |
| QL3S | 8360 | 149.8 | 86.0 | 501.3 | 232.0 | 651.1 | 249.1 | 612 | 147.5 | 83.5 | 238.6 | 227.4 | 86 | 78.9 | 49.8 |
| QL4S | 8418 | 147.3 | 84.5 | 499.9 | 229.1 | 647.2 | 246.8 | 640 | 149.4 | 77.4 | 221.8 | 226.0 | 102 | 73.5 | 48.3 |
| QL5S | 8419 | 148.2 | 86.2 | 509.7 | 232.5 | 657.9 | 249.5 | 660 | 151.9 | 84.8 | 230.4 | 222.2 | 104 | 68.1 | 50.9 |

The data of service level from the point of view of passenger are given in Table 7. All the control methods outperform the no control scheme. Q-learning algorithms increase the total number of passengers who have finished their trips at the end of simulation. One remarkable feature shown in the data is that Q-learning algorithms reduce the average waiting time and the related standard deviations comparing with the terminal station controls. Though to reduce the average waiting time and its variance is not the direct objective of our algorithm, the above observation is very uplifting coupled with the results given in Table 6.

**5.3.5. Impact of Action Sets and Comparison of the ESH and DCH**

In what follows, two important issues will be considered. One is the impact of different action sets on the performance of QL3S. The other is the impact of "substituting the Expected System Headway (ESH) for the Dynamic Circle Headway (DCH) in the cost of action" on the performances of various modified Q-learning algorithms with different numbers of successive stages to be looked ahead.

Table 8

The indices of stability and interference with the different action sets in QL3S.

| Action Set | $\sum_{t\in \bar{T}} \sigma_H(t)$ | $\bar{c}_H$ | $\sigma_{\bar{c}}$ | $\max_{t\in \bar{T}} \sigma_H(t)$ | $\min_{t\in \bar{T}} \sigma_H(t)$ | $a_\Sigma$ | $\mathbb{C}_a$ | $\bar{a}$ | $\sigma_{\bar{a}}$ | $n_M$ |
|---|---|---|---|---|---|---|---|---|---|---|
| $A_{2\times 5}$ | 20999.70 | 20.46 | 5.25 | 40.02 | 9.71 | 4366 | 3977 | 4.25 | 4.38 | 1026 |
| $A_{3\times 4}$ | 19718.70 | 19.37 | 4.95 | 39.83 | 8.28 | 5010 | 4592 | 4.92 | 5.16 | 1018 |
| $A_{4\times 4}$ | 18921.89 | 18.71 | 4.80 | 41.01 | 7.04 | 5864 | 5350 | 5.80 | 6.37 | 1011 |
| $A_{5\times 3}$ | 19520.45 | 19.32 | 4.89 | 40.27 | 9.60 | 5860 | 5409 | 5.80 | 6.22 | 1010 |
| $A_{5\times 4}$ | 17307.91 | 17.19 | 4.71 | 41.58 | 6.90 | 6015 | 5504 | 5.97 | 6.87 | 1007 |

Table 9

The indices related to waiting and riding times with the different action sets in QL3S.

| Action Set | $P_1$ | | | | | | | $P_2$ | | | | | $P_3$ | | |
|---|---|---|---|---|---|---|---|---|---|---|---|---|---|---|---|
| | $n_{P_1}$ | $t^W_{P_1}$ | $\sigma^W_{P_1}$ | $t^R_{P_1}$ | $\sigma^R_{P_1}$ | $t^{Tr}_{P_1}$ | $\sigma^{Tr}_{P_1}$ | $n_{P_2}$ | $t^W_{P_2}$ | $\sigma^W_{P_2}$ | $t^R_{P_2}$ | $\sigma^R_{P_2}$ | $n_{P_3}$ | $t^W_{P_3}$ | $\sigma^W_{P_3}$ |
| $A_{2\times 5}$ | 8360 | 149.8 | 86.0 | 501.3 | 232.0 | 651.1 | 249.1 | 612 | 147.5 | 83.5 | 238.6 | 227.4 | 86 | 78.9 | 49.8 |
| $A_{3\times 4}$ | 8358 | 145.3 | 84.0 | 502.9 | 228.3 | 648.3 | 245.8 | 659 | 149.1 | 79.1 | 214.6 | 211.4 | 92 | 68.4 | 48.7 |
| $A_{4\times 4}$ | 8358 | 147.3 | 87.3 | 511.9 | 230.4 | 659.2 | 246.9 | 548 | 142.9 | 86.1 | 254.1 | 238.1 | 77 | 69.7 | 44.3 |
| $A_{5\times 3}$ | 8400 | 147.7 | 86.4 | 511.7 | 232.8 | 659.4 | 251.3 | 636 | 151.9 | 82.6 | 249.9 | 237.2 | 113 | 63.9 | 44.7 |
| $A_{5\times 4}$ | 8299 | 146.7 | 88.4 | 509.1 | 229.6 | 655.8 | 246.8 | 621 | 147.3 | 82.1 | 236.4 | 218.2 | 111 | 69.8 | 52.0 |

To consider the first issue, we will use QL3S as the representative of Q-learning algorithms with multistage look-ahead. From the data in columns 2 to 6 in Table 8, we can see that the action set with more elements and a bigger value of the largest element will lead to a better performance of system than the action set with fewer elements and a smaller value of its biggest element. For example, the performance of bus line system with action set $A_{4\times 4}$ is better than the one with $A_{3\times 4}$.

The above benefit in performance is not free. The data in columns 7 to 9 in Table 8 show that a better performance needs stronger outside interference.

Corresponding to Table 8, the indices of service level from the point of view of passenger with the different action sets in QL3S are presented in Table 9. The data of different action sets are similar to each other. There is no such a series of indices with respect to an action set which dominates the other series. For example, the average waiting times for passengers who have arrived their destinations are 145.3s, 147.3s and 146.6s with respect to action sets $A_{3\times 4}$, $A_{4\times 4}$ and $A_{5\times 4}$, respectively. But the average waiting times for passengers who still in the buses at the end of observation period are 149.1s, 142.9s and 147.3s with respect to action sets $A_{3\times 4}$, $A_{4\times 4}$ and $A_{5\times 4}$, respectively. These data don't coincide with the changing order of the stability indices of these action sets shown in Table 8.

Based on the above observation, we can obtain two points about the application of the new method. The first point is that when the system stability is the focus of application, the action set with more elements and a bigger value of the largest element should be chosen if the corresponding outside interference is acceptable. The second point is that when some indices of service level from the point of view of passenger are the focus of application, we should choose the action set corresponding to the best values of the concerned indices. In this case, the system stability and the outside interference should also be acceptable.

Table 10

The indices of stability and interference with ESH used in the cost of action.

| Methods | $\sum_{t\in \bar{T}}\sigma_H(t)$ | $\bar{c}_H$ | $\sigma_{\bar{c}}$ | $\max_{t\in \bar{T}}\sigma_H(t)$ | $\min_{t\in \bar{T}}\sigma_H(t)$ | $a_\Sigma$ | $\mathbb{C}_a$ | $\bar{a}$ | $\sigma_{\bar{a}}$ | $n_M$ |
|---|---|---|---|---|---|---|---|---|---|---|
| QL1S | 27054.59 | 26.62 | 7.74 | 44.24 | 13.95 | 4822 | 4397 | 4.73 | 4.57 | 1018 |
| QL2S | 25095.15 | 24.65 | 5.76 | 42.12 | 11.35 | 4782 | 4312 | 4.71 | 4.48 | 1016 |
| QL3S | 24303.73 | 23.71 | 5.36 | 42.41 | 11.20 | 4628 | 4168 | 4.52 | 4.42 | 1025 |
| QL4S | 28446.21 | 27.86 | 6.75 | 44.82 | 12.37 | 4850 | 4389 | 4.75 | 4.55 | 1021 |
| QL5S | 30162.09 | 29.65 | 9.42 | 51.61 | 13.22 | 4882 | 4391 | 4.80 | 4.75 | 1017 |

Table 11

The indices related to waiting and riding times with ESH used in cost of action.

| Methods | $P_1$ | | | | | | | $P_2$ | | | | | $P_3$ | | |
|---|---|---|---|---|---|---|---|---|---|---|---|---|---|---|---|
| | $n_{P_1}$ | $t^W_{P_1}$ | $\sigma^W_{P_1}$ | $t^R_{P_1}$ | $\sigma^R_{P_1}$ | $t^{Tr}_{P_1}$ | $\sigma^{Tr}_{P_1}$ | $n_{P_2}$ | $t^W_{P_2}$ | $\sigma^W_{P_2}$ | $t^R_{P_2}$ | $\sigma^R_{P_2}$ | $n_{P_3}$ | $t^W_{P_3}$ | $\sigma^W_{P_3}$ |
| QL1S | 8397 | 147.6 | 86.6 | 508.5 | 227.9 | 656.1 | 244.7 | 659 | 158.9 | 98.2 | 240.5 | 229.3 | 86 | 73.8 | 58.0 |
| QL2S | 8455 | 149.2 | 85.4 | 505.7 | 227.7 | 654.9 | 244.8 | 660 | 149.1 | 83.3 | 235.2 | 216.0 | 96 | 67.7 | 52.3 |
| QL3S | 8257 | 143.5 | 83.8 | 500.4 | 226.0 | 643.9 | 241.7 | 649 | 152.8 | 83.5 | 242.8 | 216.2 | 108 | 91.3 | 58.3 |
| QL4S | 8399 | 149.1 | 87.5 | 503.5 | 224.5 | 652.5 | 242.5 | 650 | 152.2 | 81.8 | 222.0 | 218.8 | 69 | 67.2 | 42.7 |
| QL5S | 8493 | 145.9 | 84.8 | 504.2 | 230.0 | 650.2 | 246.0 | 617 | 149.9 | 90.1 | 229.3 | 218.4 | 98 | 81.7 | 63.3 |

The results about the second important issue are summed up in Tables 10 and 11. Comparing with the data in Tables 6 and 7, three points about the results using ESH is worth noting below.

The first point is that the performance of our holding method with ESH replacing DCH in the cost of action is still very effective. The performance shown in Table 10 is better than the one based on the terminal station controls. Bus bunching can be removed from the operation of bus line system thoroughly in this situation.

The second point is that with the increasing number of successive stages to be looked ahead, a similar trend is shown in the data related to ESH as in the data related to DCH. The operation performance of bus line system is improved at beginning, and then starts to decrease gradually. Note that the changing point in Table 10 is at QL3S that is the same as in Table 6.

The third point is that with the same depth of look-ahead, the performance of a holding method with DCH used to calculate the cost of action is better than the one with ESH. Based on this point, we always suggest using DCH in the cost of action when the modified Q-learning algorithm is implemented.

The reason behind the above observations is that DCH is better than ESH at capturing the dynamic stability feature of a bus line system. Since the goal of the control scheme with DCH used to calculate the cost of action is to stabilize a high frequency bus line with an average headway near to ESH, the performances with respect to these two headways are similar.

**5.4. Applying Different Methods to Different Bus Lines**

Table 12

The indices of stability and interference of different bus lines.

| Bus Line | Methods | $\sum_{t\in\bar{T}}\sigma_H(t)$ | $\bar{c}_H$ | $\sigma_{\bar{c}}$ | $\max_{t\in\bar{T}}\sigma_H(t)$ | $\min_{t\in\bar{T}}\sigma_H(t)$ | $a_{\Sigma}$ | $\bar{a}$ | $\sigma_{\bar{a}}$ | $n_M$ | Bunch |
|---|---|---|---|---|---|---|---|---|---|---|---|
| L1 | NC | 81136.35 | 200.33 | 81.29 | 333.65 | 32.00 | / | / | / | / | Yes |
|  | SP | 23598.19 | 61.77 | 19.18 | 107.8 | 17.40 | 1638 | 74.44 | 61.86 | 22 | No |
|  | QL1S | 8500.46 | 22.19 | 8.62 | 45.54 | 4.05 | 1638 | 4.27 | 4.45 | 383 | No |
|  | QL3S | 7526.86 | 19.65 | 7.26 | 39.52 | 3.15 | 1504 | 3.92 | 4.41 | 383 | No |
| L2 | NC | 145089.21 | 255.43 | 107.48 | 401.37 | 39.62 | / | / | / | / | Yes |
|  | SP | 43699.86 | 79.16 | 27.30 | 149.67 | 26.50 | 1596 | 69.39 | 67.40 | 23 | No |
|  | QL1S | 13410.41 | 24.83 | 7.52 | 46.51 | 10.49 | 2304 | 4.26 | 4.57 | 540 | No |
|  | QL3S | 11426.15 | 21.19 | 7.34 | 44.59 | 8.12 | 2158 | 4.00 | 4.50 | 539 | No |
| L3 | NC | 161196.65 | 221.72 | 86.23 | 364.83 | 38.37 | / | / | / | / | Yes |
|  | SP | 60547.09 | 85.03 | 30.94 | 145.91 | 25.90 | 2182 | 87.26 | 75.90 | 25 | Yes |
|  | QL1S | 18275.87 | 26.18 | 8.74 | 47.07 | 10.22 | 3282 | 4.70 | 4.56 | 698 | No |
|  | QL3S | 16022.96 | 23.02 | 6.08 | 42.67 | 10.02 | 3106 | 4.46 | 4.49 | 696 | No |
| L4 | NC | 243527.37 | 264.99 | 250.67 | 1941.09 | 38.43 | / | / | / | / | Yes |
|  | SP | 117779.82 | 138.89 | 74.73 | 284.40 | 37.71 | 6075 | 253.11 | 233.81 | 24 | Yes |
|  | QL1S | 25775.62 | 29.45 | 11.82 | 56.81 | 10.16 | 3914 | 4.47 | 4.47 | 875 | No |
|  | QL3S | 21376.14 | 24.45 | 6.83 | 42.96 | 9.58 | 3808 | 4.35 | 4.47 | 874 | No |

Table 13

The indices of service level of 4 bus lines under different control schemes.

| Bus line | Methods | $P_1$ | | | | | | | $P_2$ | | | | | $P_3$ | | |
|---|---|---|---|---|---|---|---|---|---|---|---|---|---|---|---|---|
|  |  | $n_{P_1}$ | $t^W_{P_1}$ | $\sigma^W_{P_1}$ | $t^R_{P_1}$ | $\sigma^R_{P_1}$ | $t^{Tr}_{P_1}$ | $\sigma^{Tr}_{P_1}$ | $n_{P_2}$ | $t^W_{P_2}$ | $\sigma^W_{P_2}$ | $t^R_{P_2}$ | $\sigma^R_{P_2}$ | $n_{P_3}$ | $t^W_{P_3}$ | $\sigma^W_{P_3}$ |
| L1 | NC | 3685 | 283.9 | 192.8 | 492.7 | 226.6 | 776.6 | 306.4 | 290 | 452.7 | 237.5 | 246.3 | 222.4 | 187 | 251.0 | 179.8 |
|  | SP | 3796 | 318.2 | 311.6 | 521.9 | 243.7 | 840.2 | 405.5 | 320 | 438.5 | 462.5 | 200.7 | 237.3 | 229 | 532.6 | 400.5 |
|  | QL1S | 3815 | 217.3 | 133.6 | 523.1 | 237.3 | 740.5 | 273.0 | 298 | 232.5 | 188.1 | 220.6 | 230.4 | 80 | 146.2 | 128.1 |
|  | QL3S | 3787 | 239.2 | 197.0 | 517.0 | 235.6 | 756.2 | 315.8 | 305 | 296.9 | 258.4 | 226.6 | 234.7 | 125 | 247.9 | 206.9 |
| L2 | NC | 5075 | 273.2 | 182.2 | 489.6 | 222.9 | 762.8 | 282.9 | 399 | 368.3 | 245.8 | 201.0 | 215.4 | 197 | 199.7 | 182.4 |
|  | SP | 5255 | 224.9 | 155.9 | 499.6 | 233.3 | 724.5 | 277.0 | 383 | 230.1 | 218.5 | 239.1 | 223.9 | 156 | 237.5 | 241.6 |
|  | QL1S | 5257 | 204.3 | 136.9 | 506.9 | 234.8 | 711.2 | 265.3 | 401 | 237.6 | 199.8 | 242.5 | 216.7 | 116 | 165.9 | 135.7 |
|  | QL3S | 5304 | 218.4 | 212.7 | 511.9 | 231.9 | 730.3 | 306.2 | 406 | 319.4 | 417.5 | 251.3 | 219.0 | 169 | 414.6 | 458.4 |
| L3 | NC | 6137 | 269.9 | 183.9 | 506.2 | 230.8 | 776.1 | 304.6 | 469 | 394.4 | 262.1 | 236.9 | 213.9 | 304 | 289.2 | 234.9 |
|  | SP | 6141 | 196.8 | 133.1 | 518.6 | 246.4 | 715.4 | 288.9 | 516 | 242.4 | 170.7 | 265.6 | 245.1 | 156 | 172.9 | 144.2 |
|  | QL1S | 6175 | 169.5 | 110.1 | 524.7 | 240.2 | 694.3 | 269.9 | 454 | 188.4 | 129.3 | 257.6 | 226.3 | 161 | 177.6 | 163.7 |
|  | QL3S | 6145 | 181.4 | 126.9 | 533.8 | 240.2 | 715.2 | 275.4 | 528 | 216.8 | 166.2 | 235.6 | 240.1 | 121 | 212.9 | 196.6 |
| L4 | NC | 6709 | 223.1 | 162.6 | 482.9 | 216.1 | 706.1 | 275.4 | 470 | 375.2 | 249.4 | 214.3 | 197.5 | 305 | 233.4 | 182.5 |
|  | SP | 6876 | 205.8 | 155.4 | 518.4 | 269.1 | 724.2 | 311.9 | 507 | 175.4 | 108.0 | 253.4 | 228.9 | 133 | 95.5 | 76.4 |
|  | QL1S | 6932 | 148.9 | 86.3 | 497.8 | 224.4 | 646.8 | 242.7 | 554 | 149.7 | 85.7 | 213.0 | 211.8 | 87 | 88.6 | 61.8 |
|  | QL3S | 6905 | 148.3 | 83.5 | 505.1 | 231.3 | 653.5 | 246.5 | 544 | 149.2 | 83.7 | 244.0 | 221.1 | 75 | 70.3 | 57.1 |

In the preceding subsection, we have investigated the results related to bus line L5. To learn about the effects of various holding methods on different bus lines, the main results of four other bus lines under different holding control schemes are presented in Tables 12 and 13. The basic information of these four bus lines has been summarized in subsection 5.2. The Expected System Headways of L1, L2, L3 and L4 are 311.0s, 295.96s, 286.73s and 275.66s, respectively.

Three points are worth noting here. The first point is as follows. For a given bus line, the data in Tables 12 and 13 show a similar changing trend of performance indices as in Tables 6 and 7. The effect of QL3S is better than the one of QL1S. The effect of QL1S is better than the one of the single terminal station control. They all outperform the no control scheme. To save the space, we omit here the results of some methods including two points control, QL2S, QL4S and QL5S. In fact, they have the similar performances as in Tables 6 and 7.

The second point is as follows. The average waiting times shown in Table 13 for QL1S and QL3S do not always shown the same changing trend as the average instability indices shown in Table 12. The reason for the above mismatching phenomenon is very complex. The direct reason lies in the size of holding time and the spatial-temporal feature of the holding operation. Here the spatial-temporal feature means when and where to implement the specified holding operation. As an important factor, the different demand patterns may lead to the above mismatching results. The above observation reminds us that to choose a proper holding strategy, we need to trade off between higher stability and longer waiting time in some cases.

The third point is as follows. The indices defined in subsection 2.2 can be used to compare the operation performances of different bus lines just as to compare the operation performances of a given bus line under different control schemes. For example, we can compare the performance of L3 under QL1S with the one of L4 under QL3S using the values of various indices in Tables 12 and 13. Based on the comparison, we can safely deduce a conclusion that the operation performance of L4 under QL3S dominates the performance of L3 under QL1S with respect to the stability indices and the outside interference strength.

**5.5. Computation Times of Different Methods**

One aspect of the practicality of a holding method is its operational efficiency. Before starting a holding operation off, the computation time required to supply a decision, i.e. determining a holding time, is the key to the successful implementation of our holding method.

The running environment of our program is as follows. The program used in this section is written in Java 1.8.0_45 and runs in NetBeans IDE 8.0.2. The computer used is a laptop with processor of Intel® Core i3-3120M CPU @2.50GHz and installed memory (RAM) of 4.00GB (2.32GB usable).

Table 14

Computation times with regard to different bus lines.

| Method | L1 | | L2 | | L3 | | L4 | | L5 | |
|---|---|---|---|---|---|---|---|---|---|---|
| | **Sim** | **Dec** | **Sim** | **Dec** | **Sim** | **Dec** | **Sim** | **Dec** | **Sim** | **Dec** |
| QL1S | 0.421 | <0.001 | 0.796 | <0.001 | 1.218 | <0.001 | 1.655 | <0.001 | 2.371 | <0.001 |
| QL2S | 0.515 | <0.001 | 0.858 | <0.001 | 1.342 | <0.001 | 1.804 | <0.001 | 2.512 | <0.001 |
| QL3S | 0.655 | <0.001 | 1.014 | <0.001 | 1.685 | <0.001 | 2.574 | <0.001 | 3.588 | <0.001 |
| QL4S | 1.872 | 0.015 | 3.448 | 0.015 | 5.68 | 0.015 | 8.658 | 0.015 | 12.841 | 0.015 |
| QL5S | 10.228 | 0.031 | 19.034 | 0.043 | 33.092 | 0.053 | 50.248 | 0.062 | 72.810 | 0.078 |

In Table 14, computation times with regard to different bus lines are presented. Note that the time unit is second (s) in this table. Any column with title "Sim" holds the computation times for different methods to run one time of the simulation system. Any column with title "DEC" holds the computation times for different methods to supply a control action. If the content is <0.001, it means the computation time is less than 0.001 second. The reason why we use <0.001 but not a specified number is that the result based on the function of System.currentTimeMillis() in Java language is 0.0. This means the running time is less than one millisecond.

When the number of successive stages to be looked ahead in the Q-learning algorithm is equal to or less than 3, the computation time required to supply a decision (a specified holding time) is less than 1 millisecond. When the number of successive stages increases to 4, the computation time becomes 0.015 second. Even when the number of successive stages is improved to 5, the worst case related to the longest bus line L5 is about 0.078 second. In view of the running environment of our program and the negligible computation time for making a decision, these results guarantee the practicability of our method in the real life.

To run one time of the simulation system is to mimic the operation of bus line in two hours. The computation times required to run one time of the simulation system given in Table 14 are very promising. When the number of successive stages to be looked ahead in the Q-learning algorithm is equal to or less than 3, the computation time required to run one time of the simulation system is less than 4s. The worst case with 5 successive stages is about 72.810s.

Why are the computation times required to supply a holding time in the row of QL5S different? The reason is that a long bus line with more components will occupy more computer memory than a short bus line with fewer components. This will influence the computation speed.

## 6. Discussion and Conclusions

The algorithm proposed in this study is in fact a general framework which can be used to create many specific methods. These methods are characterized by their specified number of successive stages to be looked ahead and the specified headway used to estimate the cost of action. The components of a specific method required to be determined carefully include the action set, the number of successive stages to be looked ahead, the initial value of epsilon used in the Epsilon-greedy rule, and the specific headway used to estimate the cost of action. Each choice has its own pros and cons. For example, choosing an action set with many elements may improve the operation performance of bus line system at the cost of increasing the computation time of supplying a holding time. A practical method should be determined after carefully balancing all the pros and cons.

One remarkable feature of the Q-learning algorithm is that the on-line and off-line versions of this algorithm can be integrated. To do so will strengthen the flexibility of this algorithm. In our case, the off-line version will use a simulation system as the basis of state transition to simulate the evolution of a bus line system. But the on-line version will take the actual state transition during the system evolution. Just like the off-line version, the on-line version will use some expected values to roll the simulation system forward so as to determine the optimal action. One reasonable application of these two versions could be as follows. When the field data are abundant at the beginning, the off-line version can be implemented first to get a good approximate Q-factor and then the on-line version can be adopted to further refine the approximation. But if the field data are limited at the beginning, we should adopt the on-line version directly. In this case, we had better limit the available actions to a relative small range at beginning so as to avoid the possible strong instability due to the unexpected outside interference.

The main contents of this paper can be summed up as follows. In order to cope with the complex real situation and dynamically adjust the control actions with the evolution of a bus line system, a holding strategy based on the approximate dynamic programming is proposed. A modified Q-learning algorithm with multi-stage look-ahead is designed to carry out the above holding strategy. Artificial neural network is used to approximate the unavailable Q-factor. The multi-stage look-ahead is adopted to get through the early phase labeled by considerable oscillation easily and quickly. Two indices of system stability are defined. One is the average of the pseudo standard deviation of instantaneous headways over all stages. The other is the deviation of the pseudo standard deviation of instantaneous headways over all stages. The stability indices of different bus lines can be compared with each other to evaluate their performances. The results of numerical experiments show that with smaller outside interference (or in another word shorter holding time) the modified Q-learning algorithms outperform the terminal station holding methods with respect to the stability indices and the average waiting times of passengers.

Except the already discussed issues, many other issues can be investigated in a similar way in the future. For example, holding strategy in a multiline transit system may be realized with ADP. The other interesting direction is using ADP to study the possibility and effectiveness of the cruising speed control.

**Acknowledgment**

The authors are grateful for the valuable comments of the two anonymous reviewers, the Associate Editor and the Editor in Chief Prof. Martin Savelsbergh of Transportation Science. This research was supported in part by National Natural Science Foundation of China(71601118, 71801153, 71871144），the Natural Science Foundation of Shanghai(18ZR1426200), the Shanghai First-class Academic Discipline Project (S1201YLXK) and the Key Climbing Project of USST (SK17PA02).